# Hessian measures II

By Neil S. Trudinger* and Xu-Jia Wang

## Abstract

In our previous paper on this topic, we introduced the notion of $k$-Hessian measure associated with a continuous $k$-convex function in a domain $\Omega$ in Euclidean $n$-space, $k = 1, \cdots, n$, and proved a weak continuity result with respect to local uniform convergence. In this paper, we consider $k$-convex functions, not necessarily continuous, and prove the weak continuity of the associated $k$-Hessian measure with respect to convergence in measure. The proof depends upon local integral estimates for the gradients of $k$-convex functions.

## 1. Introduction

In the paper [25], we introduced the notion of $k$-Hessian measure as a Borel measure associated to certain continuous functions, (called $k$-convex), in Euclidean $n$-space, $\mathbf{R}^n$, through the action of the $k$-Hessian operator $F_k$, $k = 1, \cdots, n$. Our results extended the special case, $k = n$, of Monge-Ampère measures associated with convex functions [1], [2], [7]. In this paper, we treat the more general setting of upper semi-continuous functions, thereby bringing our results into line with the special case, $k = 1$, of subharmonic functions in classical potential theory. The basic result in [25] was the weak continuity of the Hessian measures with respect to local uniform convergence. In this paper we prove a stronger result, (when $k \leq n/2$), namely the weak continuity of the Hessian measures with respect to local $L^1$ convergence. Our proof rests upon integral estimates, substantially different from those in [25], and we were guided somewhat in our investigations by some aspects of the theory of plurisubharmonic functions in several complex variables (see, for example, [3], [4], [6], [14], [18]). However, the analogous weak continuity result (which

---

*Research of the first author supported by Australian Research Council Grant and Humboldt Award.

AMS 1991 Mathematics Subject Classification: 58C35, 28A33, 35J60, 31B15.

would entail the weak continuity of the complex Monge-Ampère operator with respect to $L^1$ convergence) is not valid in the plurisubharmonic case and our key estimates would not be applicable there.

We shall adopt definitions and terminology similar to those introduced in [25]. For $k = 1, \cdots, n$ and $u \in C^2(\Omega)$ the $k$-Hessian operator, $F_k$, is defined by

$$(1.1) \qquad F_k[u] = S_k(\lambda(D^2u)),$$

where $\lambda = (\lambda_1, \cdots, \lambda_n)$ denotes the eigenvalues of the Hessian matrix of second derivatives, $D^2u$, and $S_k$ is the $k^{\text{th}}$ elementary symmetric function on $\mathbf{R}^n$, given by

$$(1.2) \qquad S_k(\lambda) = \sum_{i_1 < \cdots < i_k} \lambda_{i_1} \cdots \lambda_{i_k}.$$

Alternatively, we may write

$$(1.3) \qquad F_k[u] = [D^2u]_k,$$

where for any $n \times n$ matrix $\mathcal{A}$, $[\mathcal{A}]_k$ denotes the sum of its $k \times k$ principal minors. A $k$-convex function is a function which is subharmonic with respect to the operator $F_k$. A precise definition can be made in various equivalent ways. For the purpose of this introduction, we adopt a "viscosity" definition ([10], [20]). Namely, an upper semi-continuous function, $u : \Omega \to [-\infty, \infty)$, is called $k$-*convex* in $\Omega$ if $F_k[q] \geq 0$ for all quadratic polynomials $q$ for which the difference $u - q$ has a finite local maximum in $\Omega$. We will also call a $k$-convex function *proper* if it does not assume the value $-\infty$ identically on any component of $\Omega$ and denote the class of proper $k$-convex functions in $\Omega$ by $\Phi^k(\Omega)$. Note that in [25], we used the notation $\Phi^k(\Omega)$ for the subclass of *continuous* $k$-convex functions in $\Omega$. When $k = 1$, the above definition is equivalent to the usual definition of subharmonic function, with $F_1[u] = \Delta u$ for $u \in C^2(\Omega)$ (see, for example, [14]). Also $\Phi^k(\Omega) \subset \Phi^j(\Omega)$ for $j \leq k$, and a function $u \in \Phi^n(\Omega)$ if and only if it is convex on each component of $\Omega$. Furthermore, a function $u \in C^2(\Omega)$ is $k$-convex in $\Omega$ if and only if the differential inequalities

$$(1.4) \qquad F_j[u] \geq 0, \quad j = 1, \cdots, k,$$

hold in $\Omega$. This latter characterization was the basis for our definition by approximation in [25] and will be amplified, along with other properties of $k$-convex functions, in the next section. In particular, any $k$-convex function in $\Omega$ is the pointwise limit of a decreasing sequence of functions in $\Phi^k \cap C^2(\Omega')$ for any $\Omega' \subset\subset \Omega$ (Lemma 2.4).

Corresponding to Theorem 1.1 in [25], we shall establish the following characterization of $k$-Hessian measures on $\Phi^k(\Omega)$.

THEOREM 1.1. *For any $u \in \Phi^k(\Omega)$, there exists a Borel measure $\mu_k[u]$ in $\Omega$ such that*

(i) $\mu_k[u] = F_k[u]$ *for $u \in C^2(\Omega)$, and*
(ii) *if $\{u_m\}$ is a sequence in $\Phi^k(\Omega)$ converging locally in measure to a function $u \in \Phi^k(\Omega)$, the sequence of Borel measures $\{\mu_k[u_m]\}$ converges weakly to $\mu_k[u]$.*

Theorem 1.1 provides an approximation result which is fundamental for further development of the theory of the operator $F_k$. In particular it can be applied to boundary value problems and is the basis for development of the potential theoretic study of these operators. It can also be applied to the theory of curvature measures. Note that from well known properties of subharmonic functions, [14], [18], we have the inclusion, $\Phi^k(\Omega) \subset \Phi^1(\Omega) \subset L^1_{loc}(\Omega)$, and convergence in measure is equivalent to convergence in $L^1_{loc}(\Omega)$. Furthermore, Theorem 1.1 does not improve Theorem 1.1 of [25], when $k > n/2$, as then functions in $\Phi^k(\Omega)$ satisfy a local Hölder estimate [25] and the sequence $\{u_m\}$ will converge locally uniformly.

The plan of this paper is as follows. In the next section we establish various properties of $k$-convex functions, in particular equivalent characterisations corresponding to classical subharmonic functions (Lemma 2.1), as distributions (Lemma 2.2), and by approximation (Lemmas 2.3, 2.4). In Section 3, we establish local integral estimates for Hessian operators (Theorem 3.1), while in Section 4, we establish local $L^p$ estimates for $k$-convex functions and their gradients with respect to lower order Hessian operators (Theorems 4.1 and 4.3). The proof of the weak continuity result, Theorem 1.1, is then completed in Section 5 (Theorem 5.1), together with a more general result (Theorem 5.2), pertaining to mixed Hessian measures. Finally, in Section 6, we remark on the application to the Dirichlet problem for $k$-Hessian measures, although a full treatment, together with applications to capacity is deferred until a later work [26]. We also defer the treatment of signed measures (as introduced in [25] for the continuous case), as we shall approach them through the more general theory of mixed Hessian measures [26].

The authors wish to acknowledge the support of the Mathematics Institute, University of Tübingen, where this work was completed, and are particularly grateful to Gerhard Huisken for his encouragement and interest.

## 2. Properties of $k$-convex functions

In this section, we establish equivalent criteria for $k$-convexity, in particular, in terms of approximation by smooth functions, analogous to our definition for the continuous case [25]. As remarked in the introduction, a function $u \in C^2(\Omega)$ is $k$-convex in $\Omega$ if $F_j[u] \geq 0$ in $\Omega$ for $j = 1, \cdots k$; that is, the eigenvalues $\lambda = \lambda(D^2 u)$, of the Hessian matrix $D^2 u$, lie in the closed convex cone in $\mathbf{R}^n$ given by

$$(2.1) \qquad \Gamma_k = \{\lambda \in \mathbf{R}^n \mid S_j(\lambda) \geq 0, \quad j = 1, \cdots, k\}.$$

To see this, we first observe that, if $u$ is $k$-convex in $\Omega$, we must have

$$(2.2) \qquad [D^2 u + \eta]_k \geq 0$$

for any nonnegative matrix $\eta \in \mathbf{S}^{n \times n}$, whence $S_k(\lambda + \eta) \geq 0$, for any $\eta \in \mathbf{R}^n$, $\eta_i \geq 0$, $i = 1, \cdots n$. By means of the expansion,

$$S_k(\lambda_1 + \eta_1, \lambda_2, \cdots, \lambda_n) = S_k(\lambda) + \eta_1 S_{k-1}(0, \lambda_2, \cdots, \lambda_n),$$

we then infer

$$(2.3) \qquad S_{k-1;i}(\lambda) := D_i S_k(\lambda) = S_{k-1}(\lambda)\big|_{\lambda_i = 0} \geq 0,$$

for $i = 1, \cdots, n$, and consequently

$$(2.4) \qquad S_{k-1}(\lambda) = \frac{1}{n-k+1} \sum_{i=1}^{n} S_{k-1;i}(\lambda) \geq 0.$$

By replacing $\lambda$ by $\lambda + \eta$, $\eta_i \geq 0$, we subsequently conclude $S_j(\lambda) \geq 0$, $j = 1, \cdots, k$.

The reverse implication follows from the basic properties of the elementary symmetric functions $S_k$ and their associated cones $\Gamma_k$ (see for example [5], [16]). In particular we note here the following alternative characterizations of the cone $\Gamma_k$:

$$
\begin{aligned}
(2.5) \qquad \Gamma_k &= \{\lambda \in \mathbf{R}^n \mid S_j(\lambda) \geq 0, j = 1, \cdots, k\} \\
&= \{\lambda \in \mathbf{R}^n \mid 0 \leq S_k(\lambda) \leq S_k(\lambda + \eta) \text{ for all } \eta_i \geq 0, \in \mathbf{R}\} \\
&= \{\lambda \in \mathbf{R}^n \mid S_k(\lambda + \eta) \geq 0, \text{ for all } \eta_i \geq 0, \in \mathbf{R}\}.
\end{aligned}
$$

The cone $\Gamma_k$ may also be equivalently defined as the closure of the component of the positivity set of $S_k$ containing the positive cone $\Gamma^+ = \{\eta \in \mathbf{R}^n \mid \eta_i > 0, i = 1, \cdots, n\}$, as is done in [5] where the convexity of $\Gamma_k$ and the concavity in $\Gamma_k$ of the function $S_k^{1/k}$ are also treated.

The above argument also shows that the operator $F_k$ is *degenerate elliptic* with respect to $k$-convex functions $u \in C^2(\Omega)$; that is, the matrix in $\mathbf{S}^{n \times n}$ given by

$$(2.6) \qquad F_k^{ij}[u] = \frac{\partial S_k}{\partial r_{ij}}(\lambda(D^2 u))$$

is nonnegative in $\Omega$, with eigenvalues $S_{k-1;i}(\lambda)$. Also, in our definition of $k$-convex function, the condition $F_k[q] \geq 0$ can be replaced by $D^2 q \in \Gamma_k$ so that $\Phi^j(\Omega) \subset \Phi^k(\Omega)$ if $j \geq k$ and in particular $k$-convex functions will be subharmonic. Moreover our definition is related to the usual definition of subharmonic functions through the following lemma.

LEMMA 2.1. *An upper semi-continuous function* $u : \Omega \to [-\infty, \infty)$ *is $k$-convex if and only if for every subdomain $\Omega' \subset\subset \Omega$ and every function $v \in C^2(\Omega') \cap C^0(\overline{\Omega})$ satisfying $F_k[v] \leq 0$ in $\Omega'$, the following implication is true*:

$$(2.7) \qquad\qquad u \leq v \ \text{ on } \ \partial\Omega' \Longrightarrow u \leq v \ \text{ in } \ \Omega'.$$

*Proof.* Suppose that $u$ is $k$-convex in $\Omega$ and the above implication is not true. Then the function $u - v$ must assume a positive maximum at some point $y \in \Omega'$ and so also does the function $u - \widetilde{v}$, where $\widetilde{v}$ is given by

$$\widetilde{v}(x) = v(x) + \varepsilon(\rho^2 - |x - y|^2)$$

for sufficiently small positive constants $\varepsilon, \rho$. Accordingly, we have

$$\begin{aligned} F_k[v](y) &= [D^2\widetilde{v}(y) + 2\varepsilon I]_k \\ &\geq (2\varepsilon)^k, \end{aligned}$$

which is a contradiction.

On the other hand, suppose that $u - q$ assumes a local maximum at $y \in \Omega$ for a quadratic polynomial $q$ satisfying $F_k[q] < 0$. Without loss of generality, we can assume the maximum is strict and, by vertical translation, that $u(y) > q(y)$ and $u \leq q$ on the boundary $\partial N$ of some neighbourhood $N$ of $y$, so that the implication (2.7) is violated. Thus $F_k[q] \geq 0$. $\qquad\square$

When $k = 1$, the differential inequality, $F_1[v] = \Delta v \leq 0$ in $\Omega'$, in Lemma 2.1, can be replaced by Laplace's equation, $\Delta v = 0$ in $\Omega'$. This can also be done for $k > 1$ provided we relinquish the smoothness requirement, $v \in C^2(\Omega')$, with the appropriate notion of weak solution (see the remark below).

It is well-known that a distribution is equivalent to a subharmonic function if and only if its Laplacian is nonnegative. To get a corresponding criterion for $k$-convexity we introduce the dual cones,

$$(2.8) \qquad\qquad \Gamma_k^* = \{\lambda \in \mathbf{R}^n \mid \langle \lambda, \mu \rangle \geq 0 \text{ for all } \mu \in \Gamma_k\},$$

which are also closed convex cones in $\mathbf{R}^n$. Note that $\Gamma_j^* \subset \Gamma_k^*$ for $j \leq k$ with

$$\Gamma_n^* = \Gamma_n = \{\lambda \in \mathbf{R}^n \mid \lambda_i \geq 0, \ i = 1, \cdots, n\}$$

and $\Gamma_1^*$ is the ray given by

$$\Gamma_1^* = \{t(1, \cdots, 1) \mid t \geq 0\}.$$

Lemma 2.2.  *A distribution $T$ on $\Omega$ is equivalent to a $k$-convex function in $\Omega$, if and only if*

$$(2.9) \qquad\qquad T(a^{ij}D_{ij}v) \geq 0$$

*for all $v \geq 0, v \in C_0^\infty(\Omega)$ and for all constant symmetric matrices $\mathcal{A} = [a^{ij}]$ with eigenvalues $\lambda(\mathcal{A}) \in \Gamma_k^*$.*

The assertion of Lemma 2.2 is equivalent to the distributions $\sum a^{ij}D_{ij}T$ being Borel measures for all constant matrices $\mathcal{A} \in \mathbf{S}^{n \times n}$ with eigenvalues $\lambda \in \Gamma_k^*$. The proof follows readily by coordinate transformations, with the positive matrix $\mathcal{A}$ being transformed into the identity, since a function will be $k$-convex if and only if it is subharmonic with respect to all operators, $L = \mathcal{A} \cdot D^2u$, $\lambda(\mathcal{A}) \in \Gamma_k^*$. A similar argument yields a further characterization of $k$-convex functions through mean value inequalities with respect to families of concentric ellipsoids. By judicious choice of the matrix $\mathcal{A}$ in (2.9), we see that the second derivatives of a $k$-convex function will be signed Borel measures for $k \geq 2$. Since proper subharmonic functions are locally integrable so also are proper $k$-convex functions and also the process of mollification can be applied to them. In particular for a spherically symmetric mollifier $\rho \in C_0^\infty(\mathbf{R}^n)$ satisfying $\rho(x) > 0$ for $|x| < 1$, $\rho(x) = 0$ for $|x| \geq 1$ and $\int \rho = 1$, the mollification, $u_h$, defined by

$$(2.10) \qquad\qquad u_h(x) = h^{-n} \int \rho(\frac{x-y}{h})u(y)dy,$$

for $0 < h < \text{dist}(x, \partial\Omega)$, has the following properties:

Lemma 2.3.  *Let $u \in \Phi^k(\Omega)$. Then $u_h \in C^\alpha(\Omega') \cap \Phi^k(\Omega')$ for any $\Omega' \subset \Omega$ satisfying $\text{dist}(\Omega', \partial\Omega) \geq h$. Moreover, as $h \searrow 0$, the sequence $u_h \searrow u$.*

*Proof.*  The $k$-convexity of $u_h$ in $\Omega'$ is an immediate consequence of Lemma 2.2. The remainder of Lemma 2.3 follows from the basic properties of mollification and subharmonic functions. $\qquad\square$

Lemma 2.3 yields a further criterion for $k$-convex functions, which are clearly preserved by decreasing sequences.

Lemma 2.4.  *A function $u : \Omega \to [-\infty, \infty)$ is $k$-convex in $\Omega$ if and only if its restriction to any subdomain $\Omega' \subset\subset \Omega$ is the limit of a monotone decreasing sequence in $C^2(\Omega') \cap \Phi^k(\Omega')$.*

From Lemma 2.3 (or 2.4) follows an extension of Lemma 2.4 in [25].

Lemma 2.5.  *Let $u_1, \cdots, u_m \in \Phi^k(\Omega)$ and $f$ be a convex, nondecreasing function in $\mathbf{R}^n$. Then the composite function $w = f(u_1, \cdots, u_m) \in \Phi^k(\Omega)$ also.*

As a further consequence of Lemma 2.3, we prove that proper $k$-convex functions are continuous for $k > n/2$. Following [11], we first define, for $0 < \alpha \leq 1$, $\sigma \geq 0$, the weighted interior norms and semi-norms on $C^0(\Omega)$,

$$(2.11) \qquad |u|_{0;\Omega}^{(\sigma)} = \sup_{x \in \Omega} d_x^\sigma |u(x)|,$$

$$[u]_{0,\alpha;\Omega}^{(\sigma)} = \sup_{x,y \in \Omega, x \neq y} d_{x,y}^{\sigma+\alpha} \frac{|u(x) - u(y)|}{|x - y|^\alpha},$$

$$|u|_{0,\alpha;\Omega}^{(\sigma)} = [u]_{0,\alpha;\Omega}^{(\sigma)} + |u|_{0;\Omega}^{(\sigma)},$$

where $d_x = \text{dist}(x, \partial\Omega)$, $d_{x,y} = \min\{d_x, d_y\}$. The following interpolation inequality is readily demonstrated.

LEMMA 2.6.   *For any $\varepsilon > 0$, $u \in C^0(\Omega) \cap L^1(\Omega)$,*

$$(2.12) \qquad |u|_{0;\Omega}^{(n)} \leq \varepsilon^\alpha [u]_{0,\alpha;\Omega}^{(n)} + C\varepsilon^{-n} \int_\Omega |u|,$$

*for some constant $C$ depending on $n$.*

We then have the following Hölder estimate for the cases $k > n/2$.

THEOREM 2.7.   *For $k > n/2$, $\Phi^k(\Omega) \subset C^{0,\alpha}(\Omega)$ for $\alpha = 2 - n/k$ and for any $\Omega' \subset\subset \Omega$, $u \in \Phi^k(\Omega)$,*

$$(2.13) \qquad |u|_{0,\alpha;\Omega'}^{(n)} \leq C \int_{\Omega'} |u|.$$

*where $C$ depends on $k$ and $n$.*

*Proof.* First let us assume $u \in \Phi^k(\Omega) \cap C^2(\Omega)$. For completeness, we repeat our argument in [25]. Fixing a ball $B = B_R(y) \subset \Omega$, we have by calculation that the function $w$ given by

$$(2.14) \qquad w(x) = C|x - y|^{2-n/k}, \quad x \neq y,$$

for constant $C$, satisfies

$$(2.15) \qquad F_k[w] = 0$$

in $\mathbf{R}^n - \{y\}$ (for all $k = 1, \cdots, n$). Consequently, from the classical comparison principle (or Lemma 2.1) in the punctured ball $B - \{y\}$, we infer the estimate

$$(2.16) \qquad u(x) - u(y) \leq \text{osc}_B u \left( \frac{|x - y|}{R} \right)^{2-n/k},$$

provided $k > n/2$, and hence for any $\sigma \geq 0$ we obtain

$$(2.17) \qquad [u]_{0,\alpha;\Omega}^{(\sigma)} \leq |u|_{0;\Omega}^{(\sigma)}.$$

The estimate (2.13) now follows by the interpolation inequality (2.12) and we conclude the full strength of Theorem 2.7 by approximation using Lemma 2.3. $\square$

*Examples.* The functions $w$ in (2.14) yield important examples of non-smooth $k$-convex functions. Indeed if we define $w_k$ by

$$(2.18) \qquad w_k(x) = \begin{cases} |x-y|^{2-n/k}, & k > n/2, \\ \log|x-y|, & k = n/2, x \neq y, \\ -|x-y|^{2-n/k}, & k < n/2, x \neq y, \end{cases}$$

with $w_k(y) = -\infty$, $k \leq n/2$, then $w_k$ is readily seen to be $k$-convex in any domain $\Omega$, with $F_k[w_k] = 0$ in $\Omega \backslash \{y\}$. These examples also show that the Hölder exponent in Theorem 2.7 cannot be improved and furnish useful guides towards local behaviour in the cases $k \leq n/2$.

*Remark.* Our definition of $k$-convex functions coincides with the notion of the inequality "$F_k[u] \geq 0$" holding in $\Omega$ in the *viscosity* sense (see [10], [20], [27]). The proof of Lemmas 2.3 and 2.4 could have been effected by employing a basic technique from viscosity theory, namely approximation by the sup-convolution, given, for $0 < \varepsilon < \text{dist}(x, \partial\Omega)$, by

$$(2.19) \qquad u_\varepsilon^+(x) = \sup_{y \in \Omega} \left( u(y) - \frac{|x-y|^2}{2\varepsilon} \right).$$

The function $u_\varepsilon^+(x)$ will be both $k$-convex and semi-convex, for $u \in \Phi^k(\Omega)$, and consequently twice differentiable almost everywhere with $F_k[u] \geq 0$, whenever the second differential exists. By mollification and the concavity of $S_k^{1/k}$ on $\Gamma_k$, we may again arrive at Lemma 2.3. More generally for $\psi \in C^0(\Omega)$, an upper semi-continuous function $u : \Omega \to [-\infty, \infty)$ satisfies "$F_k[u] \geq \psi$" in $\Omega$ in the *viscosity* sense if for all $\varphi \in C^2(\Omega)$ and local maximum points $y \in \Omega$ of $u - \varphi$, we have $F_k[\varphi](y) \geq \psi(y)$. A lower semi-continuous function $u : \Omega \to (-\infty, \infty]$ satisfies the opposite inequality "$F_k[u] \leq \psi$" in the *viscosity* sense if for all $\varphi \in C^2(\Omega) \cap \Phi^k(\Omega)$ and local minimum points $y \in \Omega$ of $u - \varphi$, we have $F_k[\varphi](y) \leq \psi(y)$. A function $u \in C^0(\Omega)$ is then a viscosity solution of the equation, "$F_k[u] = \psi$", in $\Omega$ if both $F_k[u] \geq \psi$ and $F_k[u] \leq \psi$ in the viscosity sense. This definition coincides with those in [23] and [25] restricted to continuous $\psi$. In particular it is equivalent to the equation, $\mu_k[u] = \psi$, where $\mu_k$ is the $k$-Hessian measure of $u$ as defined in [25]. Furthermore, it follows that the function $v$ in Lemma 2.1 can be replaced by a *$k$-harmonic* function in $C^0(\overline{\Omega'})$, that is, a solution $v \in C^0(\overline{\Omega'})$ of the homogeneous equation, $F_k[u] = 0$, in $\Omega'$.

## 3. Fundamental estimates

In order to approach the proof of Theorem 1.1, we need an estimate guaranteeing the local boundedness of the sequence of measures $\mu_k[u_m]$. The following theorem is the appropriate extension of Lemma 2.3 in [25].

THEOREM 3.1.   *Let $u \in \Phi^k(\Omega) \cap C^2(\Omega)$ satisfy $u \leq 0$ in $\Omega$.   Then for any subdomain $\Omega' \subset\subset \Omega$,*

$$(3.1) \qquad \int_{\Omega'} F_k[u] \leq C \left( \int_\Omega (-u) \right)^k,$$

*where $C$ is a constant depending on $\Omega$ and $\Omega'$.*

Note that since our considerations here are local and upper semi-continuous functions are locally bounded from above, there is no loss of generality in assuming $u \leq 0$ in $\Omega$ and $u \in L^1(\Omega)$. The proof of Theorem 3.1 depends on the *classical* existence theorem of Caffarelli, Nirenberg and Spruck [5] and the interior gradient bound [8], [23], which for convenience we state here.

THEOREM 3.2 ([5]).   *Let $\Omega$ be a bounded, uniformly $(k-1)$-convex domain in $\mathbf{R}^n$ with boundary $\partial\Omega \in C^\infty$ and $\varphi, \psi$ be functions in $C^\infty(\overline\Omega)$ with $\inf_\Omega \psi > 0$.  Then there exists a unique $k$-convex function $u \in C^\infty(\overline\Omega)$ solving the Dirichlet problem*

$$(3.2) \qquad F_k[u] = \psi \quad \text{in} \ \ \Omega,$$
$$u = \varphi \quad \text{on} \ \ \partial\Omega.$$

THEOREM 3.3 ([23]).   *Let $\Omega$ be a domain in $\mathbf{R}^n$ and $u \in C^2(\Omega) \cap \Phi^k(\Omega)$ satisfy*

$$(3.3) \qquad F_k[u] = \psi_0 \quad \text{in} \ \ \Omega,$$

*for some constant $\psi_0 \geq 0$.  Then for any ball $B = B_R(y) \subset \Omega$,*

$$(3.4) \qquad |Du(y)| \leq \frac{C}{R} \big(\operatorname{osc}_B u\big),$$

*where $C$ is a constant depending on $n$.*

Using the norms (2.11) and the interpolation inequality (2.12), we can improve Theorem 3.3 as follows,

COROLLARY 3.4.   *Let $\Omega$ be a domain in $\mathbf{R}^n$ and $u \in C^2(\Omega) \cap \Phi^k(\Omega) \cap L^1(\Omega)$ satisfy equation (3.3).  Then*

$$(3.5) \qquad |u|_{0,1;\Omega}^{(n)} \leq C \int_\Omega |u|,$$

*where $C$ is a constant depending on $n$.*

*Proof of Theorem* 3.1. It is enough to consider the case of concentric balls, $\Omega = B_R(y)$, $\Omega' = B_r(y)$, $r < R$ and $u \in \Phi^k \cap C^\infty(\Omega)$. By replacement of $u(x)$ by $u(x) + \delta|x|^2/2$ for $\delta \in (0,1)$, we may also assume $F_k[u] \geq \delta^k$ in $\Omega$. Let $\eta \in C^\infty(\overline{\Omega})$ satisfy $0 \leq \eta \leq 1$, $\eta = 1$ in $B_r(y)$, $\eta = \delta^k$ for $|x - y| \geq (R + 2r)/3$ and let $\widetilde{u} \in C^\infty(\overline{\Omega})$ be the unique $k$-convex solution of the Dirichlet problem

$$(3.6) \qquad\qquad F_k[\widetilde{u}] = \eta F_k[u] \ \text{ in } \ \Omega,$$
$$\widetilde{u} = 0 \ \text{ on } \ \partial\Omega,$$

as guaranteed by Theorem 3.2. By the comparison principle, we have

$$(3.7) \qquad\qquad\qquad u \leq \widetilde{u} \leq 0$$

in $\Omega$, so that in particular,

$$(3.8) \qquad\qquad\qquad \int_\Omega |\widetilde{u}| \leq \int_\Omega |u|.$$

Let $\zeta \in C_0^\infty(\Omega)$ be a further cut-off function. Then, by integration by parts,

$$(3.9) \qquad \int_\Omega \zeta F_k[\widetilde{u}] = \frac{1}{k} \int \zeta F_k^{ij}[\widetilde{u}] D_{ij}\widetilde{u}$$
$$= \frac{1}{k} \int \widetilde{u} F_k^{ij}[\widetilde{u}] D_{ij}\zeta$$
$$\leq \frac{1}{k} \max\left(|D^2\zeta|\,|\widetilde{u}|\right) \int_{\mathrm{supp}\, D^2\zeta} F^{ii}[\widetilde{u}]$$
$$= \frac{n-k+1}{k} \max\left(|D^2\zeta|\,|\widetilde{u}|\right) \int_{\mathrm{supp}\, D^2\zeta} F_{k-1}[\widetilde{u}].$$

Choosing $\zeta = 1$ in $B_{(R+r)/2}(y)$, $\zeta = 0$ for $|x - y| \geq (5R + r)/6$, $|D^2\zeta| \leq C(R-r)^{-2}$ and using Corollary 3.4, we then arrive at the estimate

$$(3.10) \qquad\qquad \int_{\Omega'} F_k[u] \leq \frac{C}{(R-r)^{n+2}} \int_\Omega F_{k-1}[u] \int_\Omega (-u)$$

for some constant $C$ depending on $n$ and $k$. By iterating, with respect to $k$, we then obtain

$$(3.11) \qquad\qquad \int_{\Omega'} F_k[u] \leq \frac{CR^n}{(R-r)^{k(n+2)}} \left(\int_\Omega (-u)\right)^k.$$

Finally by sending $\delta \to 0$ and using a standard convergence argument we obtain (3.1).                    $\square$

Theorem 3.1 may be alternatively derived by extending the function $u$ from the smaller ball rather than the functions $F_k[u]$. We cannot avoid some loss of smoothness in this approach, which nevertheless can be overcome by mollification. Technically we can proceed in various ways, the simplest of which is to invoke an existence theorem for the homogeneous equation which follows from Theorems 3.1, 3.3 and the standard Perron process.

THEOREM 3.5. *Let $\Omega$ be a bounded domain in $\mathbf{R}^n$ which is regular for the Laplacian and $\varphi$ a function in $\Phi^k(\Omega) \cap C^0(\overline{\Omega})$. Then there exists a unique function $u \in \Phi^k(\Omega) \cap C^0(\overline{\Omega}) \cap C^{0,1}(\Omega)$ solving the Dirichlet problem,*

$$(3.12) \qquad F_k[u] = 0 \quad \text{in} \ \ \Omega,$$
$$u = \varphi \quad \text{on} \ \ \partial\Omega.$$

In accordance with the remark at the end of Section 2, the equation (3.12) may be interpreted in the viscosity sense or, more generally, in the approximation sense of [23] or operator sense of [25]. Furthermore the estimate (3.5) will be applicable to the solution $u$. Returning to the proof of Theorem 3.1, we extend the function $u \in C^2(\Omega) \cap \Phi^k(\Omega)$ from the ball $\Omega' = B_r(y)$, to the ball $\Omega = B_R(y)$, by defining $\widetilde{u}$ to be the solution of the Dirichlet problem

$$(3.13) \qquad F_k[\widetilde{u}] = 0 \quad \text{in} \ \ \Omega - \overline{\Omega'},$$
$$\widetilde{u} = 0 \quad \text{on} \ \ \partial\Omega,$$
$$\widetilde{u} = u \quad \text{on} \ \ \partial\Omega'.$$

Clearly the function

$$(3.14) \qquad \varphi = \max\{u, \ \ K(|x|^2 - R^2)\}$$

will serve as a barrier for sufficiently large $K$, and the extended function $\widetilde{u} \in \Phi^k(\Omega) \cap C^{0,1}(\overline{\Omega})$, satisfies the estimate (3.5) in the shell, $\Omega - \Omega'$. By mollifying $\widetilde{u}$ we can then proceed again through the rest of the proof of Theorem 3.1. Alternatively, we can use a result of Guan [12] to obtain an extension $\widetilde{u}$ which is smooth in $\overline{\Omega} - \Omega'$. But again mollification, or some other smoothing process, is needed to get around the lack of smoothness across $\partial\Omega'$. We shall employ such an extension in Section 5 to complete the proof of Theorem 1.1.

## 4. Gradient estimates

Our proof of Theorem 1.1 also depends upon the following local gradient estimates for $k$-convex functions, which extend the local $L^p$ estimates for the gradients of subharmonic functions when $p < n/(n-1)$.

THEOREM 4.1. *Let $u \in C^2(\Omega) \cap \Phi^k(\Omega)$, $k = 1, \cdots, n$, satisfy $u \leq 0$ in $\Omega$. Then for any subdomain $\Omega' \subset\subset \Omega$, there exist the estimates*

$$(4.1) \qquad \int_{\Omega'} |Du|^q F_l[u] \leq C \left( \int_{\Omega} |u| \right)^{q+l}$$

*for all $l = 0, \cdots, k-1$, $0 \le q < \frac{n(k-l)}{n-k}$, where $C$ is a constant depending on $\Omega, \Omega', n, k, l$ and $q$.*

When $l = 0$ in Theorem 4.1, $F_0 \equiv 1$, and we have local gradient estimates

$$(4.2) \qquad \|Du\|_{L^q(\Omega')} \le C \int_\Omega |u|$$

for $q < \frac{nk}{n-k}$, where $C$ depends on $n, k, q, \Omega$, and $\Omega'$. By Lemma 2.3, we then infer that $\Phi^k(\Omega) \subset W^{1,q}_{\mathrm{loc}}(\Omega)$; that is, that functions in $\Phi^k(\Omega)$ lie in the local Sobolev space $W^{1,q}_{\mathrm{loc}}(\Omega)$. When $k = n$, we may, of course, take $q = \infty$ in (4.2). From the Sobolev imbedding theorem [11], we then have for $q > n$, that is, for $2k > n$, a Hölder estimate as in Theorem 2.7.

To prove Theorem 4.1, we introduce a broader class of operators, namely, the $p$-$k$-Hessian operators, given for $k = 1, \cdots, n$, $p \ge 2$, $u \in C^2(\Omega)$, by

$$(4.3) \qquad F_{k,p}[u] = \big[ D(|Du|^{p-2} Du) \big]_k.$$

When $k = 1$, we obtain the well-known $p$-Laplacian operator,

$$(4.4) \qquad F_{1,p}[u] = \mathrm{div}(|Du|^{p-2} Du);$$

while using the expanded form of the $p$-Hessian,

$$(4.5) \qquad D(|Du|^{p-2} Du) = |Du|^{p-2} \Big( I + (p-2) \frac{Du \otimes Du}{|Du|^2} \Big) D^2 u,$$

we have for $k = n$, the Monge-Ampère type operator,

$$(4.6) \qquad F_{n,p}[u] = (p-1)|Du|^{n(p-2)} \det D^2 u.$$

Let us call a function $u \in C^2(\Omega)$, $p$-$k$-convex in $\Omega$ if $F_{l,p}[u] \ge 0$ for all $l = 1, \cdots, k$. We then have the following relation between $k$-convexity and $p$-$k$-convexity.

LEMMA 4.2. *Let $u \in C^2(\Omega) \cap \Phi^k(\Omega)$. Then $u$ is $p$-$l$-convex for $l = 1, \cdots, k-1$ and $p - 2 \le n(k-l)/(n-k)$.*

*Proof.* At a point $y \in \Omega$, where $Du(y) \ne 0$, we fix a coordinate system so that the $x_1$ axis is directed along the vector $Du(y)$ and the remaining axes are chosen so that the reduced Hessian $[D_{ij}u]_{i,j=2,\cdots,n}$ is diagonal. It follows then that the $p$-Hessian is given by

$$(4.7) \qquad D_i\big(|Du|^{p-2} D_j u\big) = |Du|^{p-2} \begin{cases} (p-1)D_{i1}u & \text{if } j = 1, i \ge 1, \\ D_{1j}u & \text{if } i = 1, j > 1, \\ D_{ii}u & \text{if } j = i > 1, \\ 0 & \text{otherwise.} \end{cases}$$

Hence by calculation, we obtain for $l = 1, \cdots, k - 1$ at the point $y$, setting $\widetilde{\lambda}_i = D_{ii}u(y)$, $i = 1, \cdots, n$,

(4.8)

$$|Du|^{l(2-p)} F_{l,p}[u] = (p-1)\widetilde{\lambda}_1 S_{l-1;1}(\widetilde{\lambda}) + S_{l;1}(\widetilde{\lambda}) - (p-1)\sum_{i=2}^{n} S_{l-2;1,i}(\widetilde{\lambda})(D_{i1}u)^2,$$

where $S_{k;i}(\lambda) = S_k(\lambda)\big|_{\lambda_i = 0}$ as in (2.3), and $S_{k;i,j}(\lambda) = S_k(\lambda)\big|_{\lambda_i = \lambda_j = 0}$. From the $k$-convexity of $u$, we have

$$F_k[u] = \widetilde{\lambda}_1 S_{k-1;1}(\widetilde{\lambda}) + S_{k;1}(\widetilde{\lambda}) - \sum_{i=2}^{n} S_{k-2;1,i}(\widetilde{\lambda})(D_{i1}u)^2 \geq 0$$

so that using Newton's inequality, in the form

(4.9)
$$\frac{S_{k;1}}{S_{k-1;1}} \leq \frac{l(n-k)}{k(n-l)} \frac{S_{l;1}}{S_{l-1;1}},$$

we have, for

$$p - 1 \leq \frac{k(n-l)}{l(n-k)},$$

the inequality

(4.10)

$$\frac{1}{p-1}|Du|^{l(2-p)} F_{l,p}[u] \geq \widetilde{\lambda}_1 S_{l-1;1}(\widetilde{\lambda}) + \frac{S_{k;1}}{S_{k-1;1}} S_{l-1;1}(\widetilde{\lambda}) - \sum_{i=2}^{n} S_{l-2;1,i}(\widetilde{\lambda})(D_{i1}u)^2$$

$$\geq \frac{S_{l-1;1}}{S_{k-1;1}} \sum_{i=2}^{n} S_{k-2;1,i}(D_{i1}u)^2 - \sum_{i=2}^{n} S_{l-2;1,i}(D_{i1}u)^2$$

$$= \frac{1}{S_{k-1;1}} \sum_{i=2}^{n} \Big( S_{l-1;1} S_{k-2;1,i} - S_{k-1;1} S_{l-2;1,i} \Big)(D_{i1}u)^2$$

$$= \frac{1}{S_{k-1;1}} \sum_{i=2}^{n} \Big( S_{l-1;1,i} S_{k-2;1,i} - S_{k-1;1,i} S_{l-2;1,i} \Big)(D_{i1}u)^2$$

$$\geq 0,$$

again by Newton's inequality (in $n-2$ variables). Note that in the above argument, we can assume $S_{k-1;1}(\widetilde{\lambda}) > 0$ by adding, if necessary, a quadratic function to $u$. Also the proof is simpler when $l = 1$, as the terms in $D_{i1}u$, $i \neq 1$, will not be present then.                                                  □

*Proof of Theorem* 4.1. Setting

$$p^* = 1 + \frac{k(n-l)}{l(n-k)}, \quad k < n, \ l < k,$$

we obtain from Lemma 4.2 and the formula (4.8), for $2 < p < p^*$ and $u \in C^2(\Omega) \cap \Phi^k(\Omega)$,

(4.11)
$$|Du|^{l(2-p)}F_{l,p}[u] = \frac{p^* - p}{p^* - 2}F_l[u] + \frac{p-2}{p^*-2}|Du|^{l(2-p^*)}F_{l,p^*}[u]$$
$$\geq \frac{p^* - p}{p^* - 2}F_l[u],$$

and hence, for

$$q = (p-2)l < \frac{n(k-l)}{n-k},$$

we have the estimate

(4.12)
$$|Du|^q F_l[u] \leq \frac{p^*-2}{p^*-p}F_{l,p}[u].$$

Accordingly, Theorem 4.1 will follow by estimation of $F_{l,p}[u]$ in $L^1_{\text{loc}}(\Omega)$. To accomplish this, it will be convenient for us to adopt some notation from [19]. Namely, for a real $n \times n$ matrix, $\mathcal{A} = [a_{ij}]$ (not necessarily symmetric), let us write

(4.13)
$$A_k(\mathcal{A}) = [\mathcal{A}]_k,$$
$$A_k^{ij}(\mathcal{A}) = \frac{\partial}{\partial a_{ij}}[\mathcal{A}]_k.$$

Then for any vector field $g = (g_1, \cdots, g_n)$, $g_i \in C^1(\Omega)$, $i = 1, \cdots, n$, it follows that

(4.14)
$$D_i A_k^{ij}(Dg) = 0, \quad j = 1, \cdots, n,$$
$$A_k^{ij}(Dg)D_i g_j = k A_k(Dg).$$

Hence, for any nonnegative cut-off function $\eta \in C_0^2(\Omega)$, we obtain

(4.15)
$$\int_\Omega \eta F_{l,p}[u] = \int_\Omega \eta A_l(D(|Du|^{p-2}Du))$$
$$= \frac{1}{l}\int_\Omega \eta A_l^{ij}D_i(|Du|^{p-2}D_j u)$$
$$= -\frac{1}{l}\int_\Omega |Du|^{p-2}A_l^{ij}D_i \eta D_j u.$$

From (4.7), we have

(4.16)
$$A_l^{ij}D_j u = |Du|^{(l-1)(p-2)}A_l^{ij}(D^2 u)D_j u$$
$$= |Du|^{(l-1)(p-2)}F_l^{ij}[u]D_j u,$$

so that, by substituting in (4.15), we obtain

$$(4.17) \qquad \int_\Omega \eta F_{l,p}[u] = -\frac{1}{l} \int_\Omega |Du|^{l(p-2)} F_l^{ij} D_i \eta D_j u$$

$$\leq \frac{1}{l} \int_\Omega |Du|^{q+1} |D\eta| F_{l-1}[u],$$

and hence, replacing $\eta$ by $\eta^l$ and using (4.12), we obtain

$$(4.18) \qquad \int_\Omega |Du|^q \eta^l F_l[u] \leq C \max |D\eta| \int_\Omega |Du|^{q+1} \eta^{l-1} F_{l-1}[u],$$

where $C$ is the constant in (4.12). Consequently,

$$(4.19) \qquad \int_\Omega |Du|^q \eta^l F_l[u] \leq \big(C \max |D\eta|\big)^l \int_\Omega |Du|^{q+l},$$

so that the estimate (4.1) is reduced to the case $l = 0$. To handle this case, we take $l = 1$ in (4.19) with

$$q = q(1) < \frac{n(k-1)}{n-k}.$$

If $u$ is $k$-convex for $k \geq 2$, we have

$$F_2[u] = \frac{1}{2}\big((\Delta u)^2 - |D^2 u|^2\big) \geq 0$$

and hence

$$(4.20) \qquad\qquad |D^2 u| \leq \Delta u.$$

Therefore we obtain from (4.19)

$$(4.21) \qquad \int_\Omega \eta |Du|^q |D^2 u| \leq C \max |D\eta| \int_\Omega |Du|^{1+q}$$

so that

$$(4.22) \qquad \int_\Omega \eta D\big(|Du|^{1+q}\big) \leq C \max |D\eta| \int_\Omega |Du|^{1+q}$$

and thus, by the Sobolev imbedding theorem [11], and appropriate choice of $\eta$, we obtain for any subdomain $\Omega' \subset\subset \Omega$,

$$(4.23) \qquad \||Du|^{1+q}\|_{L^{n/(n-1)}(\Omega')} \leq C d_{\Omega'}^{-1} \int_\Omega |Du|^{1+q},$$

where $d_{\Omega'} = \text{dist}(\Omega', \partial\Omega)$ and, as in (4,21), (4.22), $C$ is a constant depending on $k, q$ and $n$. The estimate (4.2) now follows by interpolation or by iteration from the subharmonic case, $k = 1$. $\qquad\square$

From Theorem 4.1 we may derive corresponding estimates for the $k$-convex functions themselves.

Theorem 4.3.   *Let $u \in C^2(\Omega) \cap \Phi^k(\Omega)$, for $k \leq n/2$, satisfy $u \leq 0$ in $\Omega$. Then for any subdomain $\Omega' \subset\subset \Omega$,*

$$(4.24) \qquad \int_{\Omega'} |u|^q F_l[u] \leq C \left( \int_\Omega |u| \right)^{l+q}$$

*for all $l = 0, \cdots, k-1$, $0 \leq q < \frac{n(k-l)}{n-2k}$, where $C$ is a constant depending on $\Omega, \Omega'$, $n, k, l$ and $q$.*

*Proof.* With $\eta \geq 0, \eta \in C_0^1(\Omega)$, we estimate

$$\int_\Omega \eta^2 (-u)^q F_l[u] = \frac{q}{l} \int_\Omega \eta^2 (-u)^{q-1} F_l^{ij} D_i u D_j u - \frac{1}{l} \int_\Omega (-u)^q F_l^{ij} D_i u D_j \eta^2$$

$$\leq \frac{q(n-l+1)}{l} \int_\Omega \eta^2 (-u)^{q-1} F_{l-1} |Du|^2$$

$$+ \frac{2(n-l+1)}{l} \int_\Omega \eta (-u)^q F_{l-1} |Du| |D\eta|$$

$$\leq \frac{(q+1)(n-l+1)}{l} \int_\Omega \eta^2 (-u)^{q-1} F_{l-1} |Du|^2$$

$$+ \frac{n-l+1}{l} \int_\Omega |D\eta|^2 (-u)^{q+1} F_{l-1}.$$

Now, for any

$$p < \frac{n(k-l+1)}{n-k},$$

we have

$$\int_\Omega \eta^2 (-u)^{q-1} F_{l-1} |Du|^2 \leq \left( \int_\Omega \eta^2 F_{l-1} |Du|^p \right)^{2/p} \left( \int_\Omega \eta^2 (-u)^{\frac{p(q-1)}{p-2}} F_{l-1} \right)^{1-2/p}$$

so that if

$$q < \frac{n(k-l)}{n-2k},$$

we may choose $p$ so that

$$q^* = \frac{p(q-1)}{p-2} < \frac{n(k-l+1)}{n-2k},$$

and the estimate (4.24) follows from Theorem 4.1 by induction on $l$.   $\square$

We remark that the case $q = 1$, $l = k-1$ in Theorem 4.3, which yields a local estimate for Hessian integrals

$$(4.25) \qquad I_{k-1}[u; \Omega'] = - \int_{\Omega'} u F_{k-1}[u],$$

may also be derived from Theorem 3.1, with the aid of the extension (3.13). Taking $\Omega, \Omega'$ and $\zeta$ as in the proof of Theorem 3.1, we obtain

$$
\begin{aligned}
I_{k-1}[u; \Omega'] &\leq \int_\Omega \zeta |u| F_{k-1}[u] \\
&= \frac{1}{2(n-k+1)} \int_\Omega \zeta |u| F_k^{ij} D_{ij}(|x|^2) \\
&= -\frac{1}{2(n-k+1)} \int_\Omega |x|^2 F_k^{ij} D_{ij}(\zeta u) \\
&\leq C R^2 \int_\Omega \left\{ \zeta F_k + |D\zeta||Du|F_{k-1} + |D^2\zeta||u|F_{k-1} \right\} \\
&\leq \frac{C R^{n+k+2}}{(R-r)^{k(n+1)}} \int_\Omega (-u)^k,
\end{aligned}
$$

by virtue of (3.11) and (3.5).

## 5. Weak continuity

In this section we complete the proof of Theorem 1.1. By Lemma 2.4, any function $u \in \Phi^k(\Omega)$ is the limit in any subdomain $\Omega' \subset\subset \Omega$ of a monotone decreasing sequence $\{u_m\} \subset C^2(\Omega') \cap \Phi^k(\Omega')$. Clearly $u_m$ also converges to $u$ in $L^1(\Omega')$. The essence of the proof of Theorem 1.1 lies in the following preliminary theorem, which also serves to define $\mu_k$.

THEOREM 5.1. *Let $\{u_m\} \subset C^2(\Omega) \cap \Phi^k(\Omega)$ converge to $u \in \Phi^k(\Omega)$ in $L^1_{\mathrm{loc}}(\Omega)$. Then the sequence $\{F_k[u_m]\}$ converges weakly to a Borel measure $\mu$ in $\Omega$.*

*Proof.* Because the sequence $\{u_m\}$ is subharmonic, we can assume without loss of generality that $u_m \leq 0$ in $\Omega$. Let us fix concentric balls $B_r = B_r(y) \subset B_R(y) = B_R$ as in the proof of Theorem 3.1. The corresponding $\{\tilde{u}_m\}$, as defined by (3.13), will then converge in $L^1(B_R)$ to a function $\tilde{u} \in \Phi^k(B_R)$ which coincides with $u$ in $B_r$ and is given, in $B_R - B_r$, by

$$
(5.1) \qquad \tilde{u} = \sup\{v \in \Phi^k(B_R - \overline{B_r}) \mid v \leq 0 \text{ on } \partial B_R, \ v \leq u \text{ on } \partial B_r\}.
$$

The inequalities $v \leq 0, (u)$ on $\partial B_R, (\partial B_r)$ respectively are to be understood as

$$
\limsup_{x \to y \in \partial B_R(\partial B_r)} v(x) \leq 0, \ (u(y)),
$$

respectively. Moreover $\tilde{u} \in C^{0,1}(\overline{B_R} - \overline{B_r})$ and satisfies the equation $F_k[\tilde{u}] = 0$ in $B_R - \overline{B_r}$ together with the estimate (3.5). For $0 < h < h_0 < R - r$, let us define the mollifications $v_m = (\tilde{u}_m)_h$, $v = (\tilde{u})_h$ so that $\{v_m\} \subset \Phi^k(B_{R-h_0}) \cap C^\infty(B_{R-h_0})$ converges to $v$ in $L^1(B_\rho)$ for any $\rho \leq R - h_0$ uniformly with respect to $h$. We shall prove that the sequence $\{F_k[v_m]\}$ converges weakly in the sense of Borel measures, uniformly with respect to $h$. To accomplish this, we first let $0 < r < \rho < R - h_0$ and fix a function $\eta \in C_0^2(B_\rho)$. Then for $l, m = 1, 2, \cdots$, and

$$(5.2) \qquad w = w_t = tv_l + (1 - t)v_m, \quad 0 \le t \le 1,$$

we have, integrating by parts,

$$(5.3)$$
$$\int_\Omega \eta \big(F_k[v_l] - F_k[v_m]\big) = \int_0^1 dt \int_{B_R} \eta F_k^{ij}[w_t] D_{ij}(v_l - v_m)$$
$$= \int_0^1 dt \int_{B_R} F_k^{ij}[w_t] D_{ij}\eta \, (v_l - v_m)$$
$$\le (n - k + 1) \max |D^2\eta| \int_0^1 dt \int_{B_\rho} F_{k-1}[w_t]|v_l - v_m|.$$

We claim now that for any $\rho \in (r, R - h_0)$,

$$(5.4) \qquad \int_0^1 dt \int_{B_\rho} F_{k-1}[w_t]|v_l - v_m| \to 0$$

as $l, m \to \infty$, uniformly in $h \le h_0$. To prove (5.4), we fix $\varepsilon \in (0, 1)$ and $N$ so that for

$$(5.5) \qquad A_\varepsilon = \{x \in B_R \mid |v_l(x) - v_m(x)| > \varepsilon\},$$

we have $|A_\varepsilon| < \varepsilon$ if $l, m \ge N$. We then have

$$(5.6) \quad \int_0^1 \int_{B_\rho} F_{k-1}(v_l - v_m)^+ \le \int_0^1 \int_{B_\rho} F_{k-1}(v_l - v_m - \varepsilon)^+ + 2\varepsilon \int_0^1 \int_{B_\rho} F_{k-1}.$$

Since the sequence $\{u_m\}$ is bounded in $L^1(B_R)$, we have

$$(5.7) \qquad \int_{B_{R-h_0}} |w_t| \le \int_{B_R} \big(t|\widetilde{u}_l| + (1 - t)|\widetilde{u}_m|\big)$$
$$\le \int_{B_R} \big(t|u_l| + (1 - t)|u_m|\big)$$
$$\le \sup_m \int_{B_R} |u_m| \le K,$$

for some fixed constant $K$. Consequently, from Theorem 3.1, we obtain

$$(5.8) \qquad \int_{B_\rho} F_{k-1}[w_t] \le CK^{k-1}$$

for some constant $C$ depending on $\rho, R$ and $n$. To estimate the first term on the right-hand side of (5.6), we let $\zeta \ge 0, \zeta \in C_0^2(B_{\rho'})$ be a cut-off function, with $\rho < \rho' < R - h_0$ and $\zeta \equiv 1$ on $B_\rho$. Setting

$$z = (v_l - v_m - \varepsilon)^+,$$

we then have, for $k > 1$,

$$(5.9) \qquad \int_{B_\rho} z F_{k-1}[w_t] \leq \int_{B_{\rho'}} \zeta z F_{k-1}[w_t]$$

$$= \frac{1}{k-1} \int_{B_{\rho'}} \zeta z F_{k-1}^{ij}[w_t] D_{ij} w_t$$

$$= -\frac{1}{k-1} \int_{B_{\rho'}} F_{k-1}^{ij} D_i(\zeta z) D_j w_t$$

$$\leq \frac{1}{k-1} \left( \int_{A_\varepsilon \cap B_{\rho'}} F_{k-1}^{ij} D_i w_t D_j w_t \right)^{1/2}$$

$$\times \left( \int_{B_{\rho'}} F_{k-1}^{ij} D_i(\zeta z) D_j(\zeta z) \right)^{1/2}.$$

To estimate the last integral on the right-hand side of (5.9), we estimate, by Hölder's inequality,

$$(5.10)$$
$$\int_{A_\varepsilon \cap B_{\rho'}} F_{k-1}^{ij} D_i w_t D_j w_t$$

$$\leq (n-k+2) \left( \int_{B_{\rho'}} F_{k-2}[w_t] |Dw_t|^q \right)^{2/q} \left( \int_{A_\varepsilon \cap B_{\rho'}} F_{k-2}[w_t] \right)^{1-2/q}$$

for $2 < q < \frac{2n}{n-k}$. At this point we use Theorem 4.1 with $l = k-2$, to estimate

$$(5.11) \qquad \int_{B_{\rho'}} F_{k-2}[w_t] |Dw_t|^q \leq CK^{k-2+q}$$

for some constant $C$ depending on $\rho', R, q, k$, and $n$, and invoke the induction hypothesis that (5.4) is valid when $k$ is replaced by $k-1$. It then follows that for $k, l \geq N' \geq N$, for a further constant $N'$ depending on $\varepsilon, \rho'$,

$$(5.12) \qquad \int_0^1 dt \int_{A_\varepsilon \cap B_{\rho'}} F_{k-2}[w_t] < \varepsilon.$$

Combining (5.9)–(5.12) we get

$$(5.13)$$
$$\int_0^1 dt \int_{B_\rho} z F_{k-1}[w_t]$$

$$\leq CK^{(k-2+q)/q} \varepsilon^{1/2-1/q} \left( \int_0^1 dt \int_{B_{\rho'}} F_{k-1}^{ij} D_i(\zeta z) D_j(\zeta z) \right)^{1/2}.$$

To complete the proof, we estimate

$$(5.14) \qquad \int_{B_{\rho'}} F_{k-1}^{ij}[w_t] D_i(\zeta z) D_j(\zeta z)$$

$$\leq \int_{B_{\rho'}} F_{k-1}^{ij} D_i[\zeta(v_l - v_m - \varepsilon)] D_j[\zeta(v_l - v_m - \varepsilon)]$$

$$\leq 2 \int_{B_{\rho'}} F_{k-1}^{ij} D_i[(v_l - v_m)\zeta] D_j[(v_l - v_m)\zeta] + 2\varepsilon^2 \int_{B_{\rho'}} F_{k-1}^{ij} D_i\zeta D_j\zeta$$

$$\leq 2 \int_{B_{\rho'}} F_{k-1}^{ij} D_i[(v_l - v_m)\zeta] D_j[(v_l - v_m)\zeta]$$
$$\qquad\qquad + 2(n-k+2)\varepsilon^2 \max |D\zeta|^2 \int_{B_{\rho'}} F_{k-2}[w_t]$$

$$= -2 \int_{B_{\rho'}} \zeta^2(v_l - v_m) F_{k-1}^{ij} D_{ij}(v_l - v_m) + 2 \int_{B_{\rho'}} (v_l - v_m)^2 F_{k-1}^{ij} D_i\zeta D_j\zeta$$
$$\qquad\qquad + 2(n-k+2)\varepsilon^2 \max |D\zeta|^2 \int_{B_{\rho'}} F_{k-2}[w_t].$$

The first integral above can be estimated as follows,

$$(5.15) \qquad -\int_0^1 dt \int_{B_{\rho'}} \zeta^2(v_l - v_m) F_{k-1}^{ij} D_{ij}(v_l - v_m)$$

$$= -\int_{B_{\rho'}} \zeta^2(v_l - v_m)\big(F_{k-1}[v_l] - F_{k-1}[v_m]\big)$$

$$\leq -\int_{B_{\rho'}} (v_l F_{k-1}[v_l] + v_m F_{k-1}[v_m])$$

$$\leq CK^k,$$

by the special case, $q = 1$, of Theorem 4.3 and (5.7), where $C$ depends on $n, \rho'$ and $R$. To control the second integral we observe that the integrand vanishes outside the support of $D\zeta$, where the interior gradient bound (3.5) is applicable. Accordingly, we may estimate

$$(5.16) \qquad \int_{B_{\rho'}} (v_l - v_m)^2 F_{k-1}^{ij} D_i\zeta D_j\zeta$$

$$\leq (n-k+2) \max |D\zeta|^2 \int_{\mathrm{supp}\, D\zeta} F_{k-2}[w_t] |v_l - v_m|^2$$

$$\leq CK^2 \int_{B_{\rho'}} F_{k-2}[w_t],$$

where $C$ depends on $n, r, \rho, \rho'$, and $R$. Combining (5.14), (5.15) and (5.16) and using Theorem 3.1, we then obtain

$$(5.17) \qquad \int_0^1 \int_{B_{\rho'}} F_{k-1}^{ij}[w_t] D_i(\zeta z) D_j(\zeta z) \leq C K^k,$$

where $C$ is a constant depending on $n, r, \rho, \rho'$, and $R$. Finally combining (5.6), (5.13) and (5.17), and interchanging $l$ and $m$, we arrive at our goal, namely,

$$(5.18) \qquad \int_0^1 \int_{B_\rho} F_{k-1}[w_t] |v_l - v_m| \leq C \varepsilon^{1/2 - 1/q}$$

for $l, m \geq N'$, where $C$ is a constant depending on $r, \rho, \rho', R, q, n, k$ and $K$. Observing that $\rho'$ and $q$ can be fixed in terms of the other constants, we conclude the proof of (5.4). Consequently the sequence $F_k[v_m]$ converges weakly to a Borel measure $\widetilde{\mu}$ in $B_R$ and, since we can start with arbitrary $r < R$, we infer that $F_k[u_m]$ converges weakly to a Borel measure $\mu$ in $B_R$. Using the weak compactness arising from the bound in Theorem 3.1, or a partition of unity, we have proved Theorem 5.1. □

From Theorem 5.1 we may define the $k$-Hessian measure associated with a function $u \in \Phi^k(\Omega)$ by

$$(5.19) \qquad \int_\Omega \eta \mu_k[u] = \lim_{h \to 0} \int_\Omega \eta F_k[u_h],$$

for any $\eta \in C_0^0(\Omega)$. The full strength of Theorem 1.1 then follows by approximation.

By examination of the proof of Theorem 5.1 we may obtain a stronger result. First we note that if $u, v \in \Phi^k(\Omega) \cap C^2(\Omega)$, we have the estimate

$$(5.20) \qquad 0 \leq k F_{k-1}^{ij}[v] D_{ij} u \leq F_k[u + v] - F_k[v].$$

By applying (5.20), instead of integrating over $t$ in (5.15), we infer, in place of (5.4),

$$(5.21) \qquad \int_{B_\rho} F_{k-1}[w_t] |v_l - v_m| \to 0$$

as $l, m \to \infty$, uniformly for $0 \leq t \leq 1$, $0 < h \leq h_0$, and moreover the sequence (5.2) may be replaced by any sequence in $\Phi^k(\Omega) \cap C^2(\Omega)$, bounded in $L^1_{\text{loc}}(\Omega)$. For two functions $u, v \in \Phi^k(\Omega)$, we may then define the mixed Hessian measure $\mu_k[u; v] = u \, d\mu_{k-1}[v]$ by

$$(5.22) \qquad \int_\Omega \eta \, d\mu_k[u; v] = \int_\Omega \eta u \, d\mu_{k-1}[v]$$
$$= \lim_{h \to 0} \int_\Omega \eta u_h \, d\mu_{k-1}[v],$$

for any $\eta \in C_0^0(\Omega)$. Furthermore, we obtain, from (5.21), the following continuity result:

THEOREM 5.2. *If $\{u_m\}$, $\{v_m\}$ are two sequences in $\Phi^k(\Omega)$, bounded in $L_{\text{loc}}^1(\Omega)$ with $u_m \to u$ in $L_{\text{loc}}^1(\Omega)$, then*

$$(5.23) \qquad\qquad \int_\Omega \eta(u_m - u)d\mu_{k-1}[v_m] \to 0$$

*for any $\eta \in C_0^0(\Omega)$. Furthermore if also $v_m \to v$ in $L_{\text{loc}}^1(\Omega)$, the sequence of mixed measures $\{\mu_k[u_m; v_m]\}$ converges to $\mu_k[u; v]$ weakly.*

A more general theory of mixed Hessian measures is developed in [26], which enables us to avoid the extensions of Section 3 in the proof of Theorem 5.1.

*Examples.* Returning to the examples $w_k$ in (2.18) of non-smooth $k$-convex functions, we have

$$(5.24) \qquad \mu_k[w_k] = \begin{cases} \left(2 - \dfrac{n}{k}\right)\left[\binom{n}{k}\omega_n\right]^{1/k}\delta_y & \text{if } k \neq \dfrac{n}{2}, \\ \left[\binom{n}{k}\omega_n\right]^{1/k}\delta_y, & \text{if } k = \dfrac{n}{2}; \end{cases}$$

(see (3.2), (3.16) in [22]), so that $\mu_k$ provides an extension of the well-known fundamental solution for the Laplacian in the case $k = 1$.

## 6. Application to the Dirichlet problem

In our previous paper [25] we considered the Dirichlet problem,

$$(6.1) \qquad\qquad \mu_k[u] = \nu \quad \text{in } \Omega,$$
$$u = \varphi \quad \text{on } \partial\Omega,$$

where $\nu$ is a finite Borel measure, $\Omega$ is a uniformly $(k-1)$-convex domain (if $k > 1$), and $\varphi \in C^0(\overline{\Omega})$. The problem (6.1) was shown to be uniquely solvable with solution $u \in \Phi^k(\Omega) \cap C^0(\overline{\Omega})$ if the measure $\nu$ can be written as

$$(6.2) \qquad\qquad \nu = \nu_1 + \nu_2,$$

where $\nu_2$ has compact support in $\Omega$ and $\nu_1 \in L^1(\Omega)$, provided $k > n/2$. For the cases $k \leq n/2$, we again infer the solvability of (6.1) from Theorem 1.1, provided we assume $\nu_1 \in L^q(\Omega)$ for $q > n/2k$, with solution $u \in \Phi^k(\Omega)$, continuous near $\partial\Omega$. To show this we first note [22] that any solution $u \in \Phi^k(\Omega)$, continuous at $\partial\Omega$, of the Dirichlet problem (6.1), satisfies an *a priori* estimate

$$(6.3) \qquad\qquad \int_\Omega |u| \leq C d^n \{\max_{\partial\Omega} |\varphi| + d^{2-n/k}\big(\nu(\Omega)\big)^{1/k}\},$$

where $C$ is a constant depending on $k$ and $n$, and $d = \operatorname{diam}(\Omega)$ (also see proof below). Next suppose that the measure $\nu_2$ in the decomposition (6.2) is supported in a subdomain $\Omega' \subset\subset \Omega$. By [23], Theorem 4.1, we infer $u \in C^{0,\alpha}(\Omega - \overline{\Omega'})$, where $\alpha \le 2 - n/kq$ and $\alpha < 1$ if $k < n$, together with an interior estimate,

$$\text{(6.4)} \qquad \operatorname{osc}_{B_{\sigma R}} u \le C\sigma^\alpha \big\{ \operatorname{osc}_{B_R} u + R^{2-n/kq} \|\nu_1\|_{L^q(B_R)}^{1/k} \big\}$$

for any concentric ball, $B_{\sigma R} = B_{\sigma R}(y)$, $B_R = B_R(y) \subset \Omega - \overline{\Omega'}$, $0 < \sigma < 1$, where $C$ is a constant depending on $k, n, p$ and $\alpha$. From the interpolation inequality (2.12), we then deduce

$$\text{(6.5)} \qquad |u|^n_{0,\alpha;\Omega-\overline{\Omega'}} \le C \int_\Omega |u|$$
$$\le C\{\max_{\partial\Omega} |\varphi| + \big(\nu(\Omega)\big)^{1/k}\},$$

where $C$ is a constant depending on $k, n, q, \alpha$ and $\Omega$.

Let $\psi_m$ be a sequence of nonnegative functions in $C_0^\infty(\Omega')$ converging weakly as measures to $\nu_2$. From [23, Th. 1.1], there exists a sequence $\{u_m\} \subset \Phi^k(\Omega) \cap C^0(\overline{\Omega})$ of solutions of the Dirichlet problems

$$\text{(6.6)} \qquad \mu_k[u_m] = \nu_1 + \psi_m \ \text{ in } \ \Omega,$$
$$u = \varphi \ \text{ on } \ \partial\Omega.$$

From the local gradient estimates, Theorem 4.1, and the estimates (6.3), (6.5), there exists a subsequence $\{u_m\}$ converging in $L^1_{\text{loc}}(\Omega) \cap C^0(\Omega - \overline{\Omega'})$ to a function $u \in L^1(\Omega) \cap C^{0,\alpha}(\Omega - \overline{\Omega'}) \cap \Phi^k(\Omega)$. Fixing a further subdomain $\Omega''$ such that $\Omega' \subset\subset \Omega'' \subset\subset \Omega$, and $\varepsilon > 0$, we then have for sufficiently large $m, l$,

$$\text{(6.7)} \qquad |u_m - u_l| \le \varepsilon \quad \text{on} \ \ \partial\Omega''.$$

Using the comparison principle in the domain $\Omega - \overline{\Omega'}$, we extend (6.7) to all of $\Omega - \Omega''$. Consequently, by Theorem 1.1, we obtain $u \in \Phi^k(\Omega) \cap C^0(\overline{\Omega} - \Omega')$, together with $\mu_k[u] = \nu$ in $\Omega$, $u = \varphi$ on $\partial\Omega$.

The estimate (6.3) can be deduced simply as follows. By mollification, it suffices to prove it for $u \in C^\infty(\overline{\Omega})$ and $\partial\Omega$ a level set of $u$. By subtracting a constant we may then assume $u = 0$ on $\partial\Omega$. Let $w \in C^\infty(\overline{\Omega})$ be the $k$-convex solution of the Dirichlet problem,

$$\text{(6.8)} \qquad F_k[w] = 1 \quad \text{in} \ \ \Omega,$$
$$w = 0 \quad \text{on} \ \ \partial\Omega.$$

Then, by integration by parts,

(6.9)
$$
\begin{aligned}
-\int_\Omega u &= -\int_\Omega u F_k[w] \\
&= -\frac{1}{k}\int_\Omega w F_k^{ij}[w] D_{ij}u \\
&\le \frac{1}{k}\sup_\Omega(-w)\int_\Omega F_k^{ij}[w] D_{ij}u \\
&= \frac{1}{k}\sup_\Omega(-w)\int_\Omega D_i\big(F_k^{ij}[w] D_j u\big) \\
&= \frac{1}{k}\sup_\Omega(-w)\int_{\partial\Omega} F_k^{ij}[w]\gamma_i\gamma_j|Du| \\
&= \frac{1}{k}\sup_\Omega(-w)\int_{\partial\Omega} H_{k-1}[\partial\Omega]|Dw|^{k-1}|Du| \\
&\le \frac{1}{k}\sup_\Omega(-w)\left(\int_{\partial\Omega} H_{k-1}[\partial\Omega]|Dw|^k\right)^{1-1/k}\left(\int_{\partial\Omega} H_{k-1}[\partial\Omega]|Du|^k\right)^{1/k} \\
&= \sup_\Omega(-w)|\Omega|^{1-1/k}\big(\mu_k[u](\Omega)\big)^{1/k}
\end{aligned}
$$

by Reilly's integration formula, [19], and (6.3) follows. Alternatively for any ball $B = B_r(0) \supset \Omega$, and $\psi \in C^\infty(\overline{B})$, $\psi > F_k[u]\chi_\Omega$ in $\overline{B}$, we have $u \ge v$ in $\Omega$, where $v \in C^\infty(\overline{\Omega})$ is the solution of the Dirichlet problem

(6.10)
$$
\begin{aligned}
F_k[v] &= \psi \quad \text{in } B, \\
v &= 0 \quad \text{on } \partial B.
\end{aligned}
$$

But then, by integration by parts,

(6.11)
$$
\begin{aligned}
-\int_B v &= \frac{1}{2n}\int_B (R^2-|x|^2)\Delta v \\
&\le \frac{R^2}{2n}\int_{\partial B}(D_\gamma v) \\
&\le \frac{R^2}{2n}\big(\omega_n R^{n-1}\big)^{1-1/k}\left(\int_{\partial B}(D_\gamma v)^k\right)^{1/k} \\
&\le C R^{2+(n-1)(1-1/k)}\left(\int_B \psi\right)^{1/k},
\end{aligned}
$$

again by the Reilly formula, where $C$ depends on $n, k$ and hence (6.3) follows again as $\psi \to \chi_\Omega F_k[u]$.

The estimate (6.4) follows directly from (6.3) and the interior gradient bound (3.4) by a standard pertubation technique [23]. The assumption $\nu_1 \in L^q(\Omega)$, $q > n/2k$, can be replaced by a Morrey condition, namely for any ball $B \subset \Omega$,

$$(6.12) \qquad \nu_1(B) \leq A|B|^\delta$$

for positive constants $A$ and $\delta > 1 - \frac{2k}{n}$.

The uniqueness of the solutions to (6.1) is a more complex issue in the discontinuous case and is related to appropriate notions of capacity. In the special case, $\nu = \delta_y$ for some point $y \in \Omega$, we may infer that the Green functions $G_y$, obtained by solving (6.1), is unique [26]. Note that when $\nu_2 \in L^1(\Omega)$ in the decomposition (6.2), the uniqueness of solutions to (6.1) follows from the method in [23, Th. 2.2].

AUSTRALIAN NATIONAL UNIVERSITY, CANBERRA, ACT 0200
*E-mail addresses*: Neil.Trudinger@maths.anu.edu.au
X.J.Wang@maths.anu.edu.au